\documentclass[a4paper]{article}
\usepackage{amsmath}
\usepackage{amsthm}
\usepackage{amsfonts}
\usepackage{amssymb}
\usepackage{mathrsfs}
\usepackage{tikz}
\usepackage{environ}
\usepackage{elocalloc}
\usepackage{url}
\usepackage{caption}
\usepackage{CJKutf8}
\usepackage{mathtools}
\usepackage{dcpic}
\usepackage{tikz-cd}
\usepackage{stmaryrd}
\usetikzlibrary{positioning}
\usetikzlibrary{arrows,chains,positioning,scopes,quotes}
\usetikzlibrary{decorations.markings}
\usepackage[margin=0.9in]{geometry}

\title{Weak Counterexamples to $L^2$ Curvature Estimates for Minimizing Surfaces}
\author{Zhenhua Liu}
\begin{document}
	\maketitle

\newcommand{\ai}{\alpha}
\newcommand{\be}{\beta}
\newcommand{\Ga}{\Gamma}
\newcommand{\ga}{\gamma}	
\newcommand{\de}{\delta}
\newcommand{\De}{\Delta}
\newcommand{\e}{\epsilon}
\newcommand{\lam}{\lambda}
\newcommand{\Lam}{\Lamda}
\newcommand{\om}{\omega}
\newcommand{\Om}{\Omega}
\newcommand{\si}{\sigma}
\newcommand{\Si}{\Sigma}
\newcommand{\vp}{\varphi}
\newcommand{\rh}{\rho}
\newcommand{\ta}{\theta}
\newcommand{\Ta}{\Theta}
\newcommand{\W}{\mathcal{O}}
\newcommand{\ps}{\psi}

\newcommand{\mf}[1]{\mathfrak{#1}}
\newcommand{\ms}[1]{\mathscr{#1}}
\newcommand{\mb}[1]{\mathbb{#1}}
\newcommand{\cd}{\cdots}

\newcommand{\s}{\subset}
\newcommand{\es}{\varnothing}
\newcommand{\cp}{^\complement}
\newcommand{\bu}{\bigcup}
\newcommand{\ba}{\bigcap}
\newcommand{\co}{^\circ}
\newcommand{\ito}{\uparrow}
\newcommand{\dto}{\downarrow}

\newcommand{\ti}[1]{\tilde{#1}}
\newcommand{\la}{\langle}
\newcommand{\ra}{\rangle}
\newcommand{\ov}[1]{\overline{#1}}
\newcommand{\no}[1]{\left\lVert#1\right\rVert}
\DeclarePairedDelimiter{\cl}{\lceil}{\rceil}
\DeclarePairedDelimiter{\fl}{\lfloor}{\rfloor}
\DeclarePairedDelimiter{\ri}{\la}{\ra}

\newcommand{\du}{^\ast}
\newcommand{\pf}{_\#}
\newcommand{\is}{\cong}
\newcommand{\n}{\lhd}
\newcommand{\m}{^{-1}}
\newcommand{\ts}{\otimes}
\newcommand{\ip}{\cdot}
\newcommand{\op}{\oplus}
\newcommand{\xr}{\xrightarrow}
\newcommand{\xla}{\xleftarrow}
\newcommand{\xhl}{\xhookleftarrow}
\newcommand{\xhr}{\xhookrightarrow}
\newcommand{\mi}{\mathfrak{m}}
\newcommand{\wi}{\widehat}
\newcommand{\sch}{\mathcal{S}}

\newcommand{\w}{\wedge}
\newcommand{\X}{\mathfrak{X}}
\newcommand{\pd}{\partial}
\newcommand{\dx}{\dot{x}}
\newcommand{\dr}{\dot{r}}
\newcommand{\dy}{\dot{y}}
\newcommand{\dth}{\dot{theta}}
\newcommand{\pa}[2]{\frac{\pd #1}{\pd #2}}
\newcommand{\na}{\nabla}
\newcommand{\dt}[1]{\frac{d#1}{d t}\bigg|_{ t=0}}
\newcommand{\ld}{\mathcal{L}}

\newcommand{\N}{\mathbb{N}}
\newcommand{\R}{\mathbb{R}}
\newcommand{\Z}{\mathbb{Z}}
\newcommand{\Q}{\mathbb{Q}}
\newcommand{\C}{\mathbb{C}}
\newcommand{\bh}{\mathbb{H}}

\newcommand{\lix}{\lim_{x\to\infty}}
\newcommand{\li}{\lim_{n\to\infty}}
\newcommand{\infti}{\sum_{i=1}^{\infty}}
\newcommand{\inftj}{\sum_{j=1}^{\infty}}
\newcommand{\inftn}{\sum_{n=1}^{\infty}}	
\newcommand{\snz}{\sum_{n=-\infty}^{\infty}}	
\newcommand{\ie}{\int_E}
\newcommand{\ir}{\int_R}
\newcommand{\ii}{\int_0^1}
\newcommand{\sni}{\sum_{n=0}^\infty}
\newcommand{\ig}{\int_{\ga}}
\newcommand{\pj}{\mb{P}}

\newcommand{\io}{\textnormal{ i.o.}}
\newcommand{\aut}{\textnormal{Aut}}
\newcommand{\out}{\textnormal{Out}}
\newcommand{\inn}{\textnormal{Inn}}
\newcommand{\mult}{\textnormal{mult}}
\newcommand{\ord}{\textnormal{ord}}
\newcommand{\F}{\mathcal{F}}
\newcommand{\V}{\mathbf{V}}	
\newcommand{\II}{\mathbf{I}}
\newcommand{\ric}{\textnormal{Ric}}
\newcommand{\sef}{\textnormal{II}}

\newcommand{\wh}{\Rightarrow}
\newcommand{\eq}{\Leftrightarrow}

\newcommand{\eqz}{\setcounter{equation}{0}}
\newcommand{\se}{\subsection*}
\newcommand{\ho}{\textnormal{Hom}}
\newcommand{\ds}{\displaystyle}
\newcommand{\tr}{\textnormal{tr}}
\newcommand{\id}{\textnormal{id}}
\newcommand{\im}{\textnormal{im}}
\newcommand{\ev}{\textnormal{ev}}

\newcommand{\gl}{\mf{gl}}
\newcommand{\sll}{\mf{sl}}
\newcommand{\su}{\mf{su}}
\newcommand{\so}{\mf{so}}
\newcommand{\ad}{\textnormal{ad}}
\newcommand{\hm}{\mathcal{H}}
\newcommand{\ka}{\kappa}

\theoremstyle{plain}
\newtheorem{thm}{Theorem}[section]
\newtheorem{lem}[thm]{Lemma}
\newtheorem{prop}[thm]{Proposition}
\newtheorem*{cor}{Corollary}
\newtheorem*{pro}{Proposition}

\theoremstyle{definition}
\newtheorem{defn}{Definition}[section]
\newtheorem{conj}{Conjecture}[section]
\newtheorem{exmp}{Example}[section]

\theoremstyle{remark}
\newtheorem*{rem}{Remark}
\newtheorem*{note}{Note}

\section{Introduction}
We will prove the following theorem.
\begin{thm}
There exists a sequence of metrics $g_k$ on $\R^4$ with $\no{g_k-\de}_{C^4}\to 0$, as $k\to\infty$, where $\de$ is the standard Euclidean metric, and a sequence of smooth 2-d surfaces $\Si_k$ in a cylinder $B_{\frac{1}{2}}(0)\times \R^2$, minimizing area with respect to $g_k$, so that $\hm^2_{g_k}(\Si_k)$ is uniformly bounded but $\displaystyle\lim_{k\to\infty}\int_{\Si_k}\no{A}^2=\infty$.
\end{thm}
This provides a weak counterexample to a direct analog of the $L^2$ curvature estimates for stable minimal hypersurfaces established by Schoen and Simon in \cite{SS}.

The idea is that $L^2$ norm of $\no{A}$ is essentially controlled by genus, so if we can produce a sequence of minimizing 2-D surfaces with genus going to infinity yet bounded area, then we have proven the theorem. We use the ingenuous construction of \cite{DPH} and \cite{DPHM} to achieve this.
\subsection*{Remarks}
To the author's knowledge, this is the first written proof of a weak counterexamples to $L^2$ second fundamental form bounds for minimizing currents of higher codimension. However, the author believes that potential examples must have appeared before. Indeed, Professor De Lellis has pointed out to the author that the infinite genus example in \cite{DPHM} should have infinite $L^2$ second fundamental form. Details will appear later. 

For mod $2$ area-minimizing currents, it is straightforward to use \cite{NV} to prove uniform interior weak $L^2$ bounds. However, the author doesn't know whether this is sharp, as there are relatively few examples of minimizing currents in mod $2.$ The examples in this paper are not minimizing mod $2.$ (If so, they would converge to a mod $2$ minimizing current. However, they converge to a double copy of the disk, which is $0$ mod $2.$) 

The motivation to investigate this kind of curvature estimates lies in its broad applications. For examples, in codimension 1, the regularity part of recently popular Almgren-Pitts min-max theory is based on various curvature estimates for stable minimal hypersurfaces, including \cite{SS}. In higher codimension, the lack of such curvature bounds is a huge obstacle to establishing corresponding regularity theories. 

Finally, for clarity of presentation, the constants $C$ involved in this manuscript will not be numbered and thus possibly be different for each line, but they will be uniformly bounded.
\section{The geometric picture}
First let's explain the geometric picture. By Gauss-Bonnet, for any two dimensional surface $\Si$ of genus $g$ with smooth connected boundary $\ga,$ we have $\chi(\Si)=2-2g-1,$ and thus
\begin{align*}
\int_{\ga}\ka d\hm^1+\int_{\Si}Kd\hm^2=2\pi\chi(\Si)=2\pi-4\pi g,
\end{align*}where $\ka$ is the mean curvature or equivalently the geodesic curvature of $\ga$ in $\Si,$ and $K$ is the Gauss curvature of $\Si$. 

Suppose $\Si$ is embedded in $M^n.$ We will use $e_1,e_2$ do denote local frames that span the tangent space of $\Si.$ Let $A$ denote the second fundamental form of $\Si $ in $M^n$, $H$ the mean curvature of $\Si\s M,$ and $R$ the Riemannian curvature tensor of $M.$ Then by Gauss equations, we have
\begin{align*}
2K=2R(e_1,e_2,e_1,e_2)+\no{H}^2-\no{A}^2,
\end{align*}where $\no{}$ is the norm on tensors induced by the metric on $\Si.$
If we further assume the mean curvature $H$ of $\Si$ in $M^n$ is $0$, then we have
\begin{align*}
2K=&2R(e_1,e_2,e_1,e_2)-\no{A}^2.
\end{align*}
Collecting this into Gauss-Bonnet, we have
\begin{align*}
\int_{\Si}\no{A}^2d\hm^2=4\pi g-2\pi+\int_{\ga}\ka d\hm^1+2\int_{\Si}R(T_p\Si)d\hm^2(p).
\end{align*}
Thus, if we have control on $\hm^2(\Si),$ $\no{R}$, $|\ka|,$ and $\hm^1(\ga).$ Then $L^2$ norm of $\no{A}$ is dominated by the genus of $\Si.$
\section{The construction of the candidate surfaces}
Let $a_{k,j}=\frac{j}{3k},$ where $k\in \N$ and $-k\le j\le k$. Use $D$ to denote $B_{\frac{1}{2}}(0)$, the ball of radius $\frac{1}{2}$ around $0$ on $\C=\R^2.$ Consider the following sequence of polynomials $\{p_k\}$ defined on $D$,
\begin{align*}
p_k(z)=\prod_{j=-k}^k(z-a_{k,j}).
\end{align*}
We have $|p_k|\le (\frac{5}{6})^{2k+1}, $ $|p_k'|\le (2k+1)(\frac{5}{6})^{2k}$, and $|p_k''|\le (2k+1)2k(\frac{5}{6})^{2k-1}.$ Note that as $k\to\infty$, we have that $|p_k|,|p_k'|,|p_k''|$ all uniformly goes to $0.$ 

Consider the following sequence of maps $G_k$ from $D=B_{\frac{1}{2}}(0)\s \C$ to $D\times \C\s \R^4,$
\begin{align*}
G_k(z)=(z^2,p_k(z)).
\end{align*}
We claim the following
\begin{enumerate}
	\item  the current $(G_k)\pf[[D]]$ is area minimizing in $\C^2$, with boundary $\ga_k=(G_k)\pf[[\pd B_{\frac{1}{2}}(0)]]$
	\item  $(G_k)\pf[[D]]$ a smooth immersed surface with boundary,
	\item  $(G_k)\pf[[D]]$ is embedded outside of the transverse intersection points $(a_{k,j}^2,0),$
	\item we have $|z|\le \no{DG_k}_{\de}\le C_0,\no{D^2G_k}_{\de}\le C_0,$ with $C_0$ independent of $k.$
\end{enumerate} 
Claim 1 is direct consequence of Wirtinger's inequality. $z^2$ has nonzero derivative away from zero, while $p_k(z)=z\prod_{j\not=0}(z-a_{k,j})$ has nonzero derivative only at $0.$ Thus, by constant rank theorem we can verify claim 2. For claim 3, note that the immersed points are the noninjective points. If $z^2=w^2,$ then $z=\pm w.$ However, the solutions to $p_k(w)=p_k(-w)=-p_k(w)$ are precisely $a_{k,j}.$ To see that the intersection is transverse, note that $(G_k)\pf [[D]]$ can be represented by the restriction to $D\times \C$ of the complex algebraic variety\begin{align}\label{variety}
v^2=u\prod_{j=1}^k(u-a_{k,j}^2)^2,
\end{align}in the $(u,v)$-plane. (Let $u=z^2$ and $v=p_k.$ Then it's clear that equation (\ref{variety}) follows.) We can use this to show that the tangent cones to $(G_k)\pf[[D]]$ at $(a_{k,j}^2,0)$ are precisely $$v=\pm C_{k,j}(u-a_{k,j}^2),$$ where $C_{k,j}=\sqrt{\prod_{1\le l\le k,l\not=|j|}(a_{k,j}^2-a_{k,l}^2)}\not=0.$ This verifies claim 3. Claim 4 obviously follows from the estimates on $|p_k'|,|p_k''|$ we have established, and using Cauchy-Riemann equation repetitively to represent holomorphic derivative in real partial derivatives.

\section{Modifying the candidate surfaces and constructing the metrics}
Now, at every nonzero $(a_{k,j}^2,0),$ we substitute $(G_k)\pf[[D]]$ in each small balls $B_{r_{k,j}}((a_{k,j}^2,0))$ with a smooth neck $$(v+ C_{k,j}(u-a_{k,j}^2))(v- C_{k,j}(u-a_{k,j}^2))=\eta_{k,j}.$$ Then we use Proposition 4.4 in \cite{DPH} as in Section 5 of \cite{DPH} to glue each neck to $(G_k)\pf[[D]]$ in $B_{10r_{k,j}}$ to form a surface $\Si_k$ coincide with $(G_k)\pf[[D]]$ outside of $\cup_{j\not=0}B_{10r_{k,j}}(a_{k,j}^2,0)$ and coincide with $(v+ C_{k,j}(u-a_{k,j}^2))(v- C_{k,j}(u-a_{k,j}^2))=\eta_{k,j}$ inside $\cup_{j\not=0}B_{r_{k,j}}(a_{k,j}^2,0)$. By construction, $\Si_k$ is a smooth genus $2k$ surface with same boundary $\ga_k$ as $(G_k)\pf[[D]].$  Moreover, $\Si_k$ is the unique area minimizing surface bounded by $\ga_k$ in some $g_k,$ so that $g_k$ with coincide with $\de$ outside of $\cup_{j\not=0}B_{10r_{k,j}}(a_{k,j}^2,0)\backslash\cup_{j\not=0}B_{r_{k,j}}(a_{k,j}^2,0)$. By requiring $\eta_{k,j}<\frac{1}{100k}$ small enough, we can impose the conditions $\no{g_k-\de}_{C^4}<\frac{1}{k}.$

We have $\no{R_{g_k}}_{g_k}\le \frac{C_1}{k},$ $\no{DG_k}_{g_k}\le C_1\no{DG_k}_{\de},$ $\no{D^2G}_{g_k}\le C_1\no{D^2G}_{\de},$ with $C_1$ independent of $k.$ By the area formula, we immediately deduce that $\hm^2_{g_k}((G_k)\pf[[D]])\le C_2$, and $\hm^1_{g_k}(\ga_k)\le C_2,$ with $C_2$ independent of $k.$ Moreover, at every point on $\Ga$, we have $|\ka_k|_{g_k}\le|H_{\ga_k\s \R^4,g_k}|.$ However, $H_{\ga_k\s \R^4,g_k}$ is bounded by $\no{D G_k}\m\no{D^2 G_k},$ which is bounded by $|z|\m C.$ Since $\ga_k$ is always image of constant radius circles, this implies that we always have $\int_{\ga}^{}\ka_{g_k}d\hm^1_{g_k}<C,$ with $C$ independent of $k.$ By minimality of $\Si_k,$ we deduce that $\hm^2_{g_k}(\Si_k)\le \hm^2_{g_k}((G_k)\pf[[D]]).$ Collecting all these estimates, we have $\int_{\ga_k}\ka_k d\hm^1_k<C$ and $\int_{\Si}R_k(T_p\Si_k)d\hm^2_k(p)<\frac{C}{k}$ with $C$ independent of $k.$ By letting $k\to\infty,$ we deduce that
\begin{align*}
\int_{\Si_k}\no{A_k}^2_{g_k}d \hm^2_k\ge 8\pi k+C_0\to\infty
\end{align*}as $k\to\infty.$ This proves the theorem.
\section*{Acknowledgements}
I would first like to thank my advisor Professor Camillo De Lellis for his unwavering support during the last year, when I have been recovering from illness. I would also like to than him for suggesting this topic to me.

\end{document}